\documentclass[12pt]{article}
\title{On deformations of {\bf Q}-factorial symplectic varieties}
\author{Yoshinori Namikawa}
\date{ }
\usepackage{amsmath,amssymb,amsthm, amscd}
\chardef\bslash=`\\

\def\0{{\mathcal O}}

\begin{document}
\maketitle

\section{Introduction} 
This is a supplement to [Na]. 
In [Na] we have proposed a new category of 
complex symplectic varieties which admit 
certain singularities. It includes 
symplectic V-manifolds (cf. [Fu 1]) 
and O'Grady's singular 
moduli spaces of semi-stable torsion free sheaves 
on a K3 surface [O] (cf. [Na], Introduction). 
A normal compact K\"{a}hler space $Z$ is a 
symplectic variety if its regular locus $U$ 
admits a non-degenerate holomorphic closed 
2-form $\omega$ with the following property: 
for any resolution $\nu: \tilde{Z} \to Z$ of $Z$ 
with ${\nu}^{-1}(U) \cong U$, the 2-form $\omega$ 
on ${\nu}^{-1}(U)$ extends to a holomorphic 2-form 
on $\tilde{Z}$. If $Z$ can be embedded in a projective 
space, $Z$ is called a projective symplectic variety. 

Our purpose is to give a positive answer 
to the following problem posed in [Na]: 
\vspace{0.2cm}

{\bf Problem}. {\em Let $Z$ be a {\bf Q}-factorial 
projective symplectic variety with terminal singularities. 
Assume that $Z$ is smoothable by a suitable flat 
deformation. Is $Z$ then non-singular from the first ?} 
\vspace{0.2cm}

We actually prove more: 
\vspace{0.2cm}

{\bf Main Theorem}. {\em Let $Z$ be a {\bf Q}-factorial 
projective symplectic variety with terminal singularities. 
Then any flat deformation of $Z$ is locally trivial; in 
other words, it preserves all singularities on $Z$.}  
\vspace{0.2cm} 

This theorem would explain why ``local Torelli theorem" 
[Na, Theorem 8] holds for such symplectic varieties. 
\vspace{0.2cm}

{\bf Example} ([O], cf. also [K-L-S]): 
Let $S$ be a K3 surface with Picard 
number $1$. Fix an 
even integer $c$ with $c \geq 6$. Let $M_c$ be the moduli 
space of semi-stable sheaves on $S$ of rank 2 with $c_1 = 0$ 
and $c_2 = c$.  O'Grady proved in [O] that $M_c$ is a 
projective symplectic variety with terminal singularities. 
In the preprint version of [O]: math.AG/9708009, he constructed 
a projective resolution $\hat{\pi}: \hat{M}_c \to M_c$. 
Let $N_1(\hat{M}_c/M_c)$ be 
the {\bf R}-vector space of $\hat{\pi}$-numerical classes of 1-cycles 
$l$ on $\hat{M}_c$ such that $\hat{\pi}(l) = \mathrm{point}$.   
Then $$N_1(\hat{M}_c/M_c) = \mathbf{R}[\hat{\epsilon}_c] 
\oplus \mathbf{R}[\hat{\sigma}_c] \oplus \mathbf{R}[\hat{\gamma}_c]$$ 
for effective 1-cycles $\hat{\epsilon}_c$, $\hat{\sigma}_c$ and 
$\hat{\gamma}_c$ (ibid. Section 3). 
On the other hand, $\mathrm{Exc}(\hat{\pi})$ consists of 
three prime divisors $\hat{\Omega}_c$, $\hat{\Sigma}_c$ 
and $\hat{\Delta}_c$ (ibid. Section 3). 
This implies that $M_c$ is {\bf Q}-factorial. 
In fact, let $D$ be a prime divisor of $M_c$ and let $D'$ be 
its proper transform on $\hat{M}_c$. By the observation above, 
we can find $a_i \in \mathbf{Q}$ $(i = 1,2,3)$ such that 
$D'':= D' + a_1\hat{\Omega}_c + a_2\hat{\Sigma}_c + a_3\hat{\Delta}_c$ 
is $\hat{\pi}$-numerically trivial. Since $M_c$ has only rational 
singularities, there is 
a positive integer $m$ such that $mD''$ is linearly equivalent to the 
pull-back of a Cartier divisor $E$ of $M_c$. Then $mD$ is linearly 
equivalent to $E$. 
More generally, for a K3 surface or an Abelian surface, almost all 
moduli spaces of 
semi-stable sheaves  
with non-primitive Mukai vectors are 
{\bf Q}-factorial projective symplectic varieties [K-L-S].  
Applying Main Theorem to these moduli spaces, 
we conclude 
that any deformation of them are locally rigid; in particular, 
they have no smoothing via deformations.  
\vspace{0.2cm}
 
In the final section, we shall apply Main Theorem to get 
the following remarkable results. 
\vspace{0.2cm}

{\bf Corollary 1}. {\em Let $h: Y \to Z$ and 
$h': Y' \to Z$ be two {\bf Q}-factorial terminalizations 
of a projective symplectic variety $Z$; in other words, 
$h$ (resp. $h'$) is a crepant projective birational 
morphism such that $Y$ (resp. $Y'$) has only 
{\bf Q}-factorial terminal singularities. Assume that 
$Y$ is non-singular. Then $Y'$ is also non-singular.} 
\vspace{0.2cm}

{\bf Remark 1}. A {\bf Q}-factorial terminalization $Y$ 
of a projective symplectic variety $Z$ is again symplectic. 
In fact, there is a 2-form $h^*\omega$ on $h^{-1}(U)$,  
where $U := Z_{reg}$. By the definition of a symplectic 
variety, this 2-form extends to a holomorphic 2-form on 
$V := Y_{reg}$. Moreover, 
since $K_Y = h^*K_Z$, the extended 2-form is non-degenerate 
over $V$. If we take a resolution $\tilde{Y}$ of $Y$, then 
this 2-form further extends to a holomorphic 2-form on $\tilde{Y}$ 
\vspace{0.2cm}

{\bf Example} ({\bf flops}): When  
the codimensions of $\mathrm{Exc}(h)$ and 
$\mathrm{Exc}(h')$ are both larger than one, 
Corollary 1 says that {\em the smoothness is preserved 
under a flop}. A similar result has been proved by Kaledin  
under certain conditions 
when $Z$ is an affine symplectic 
singularity (cf. ``Symplectic resolutions: deformations 
snd birational maps", math.AG/0012008).  
\vspace{0.2cm}

{\bf Corollary 2}. {\em Let $Z$ be a projective 
symplectic variety. Assume that $Z$ has a {\bf Q}-factorial 
terminalization 
$h: Y \to Z$. Then the following are equivalent. 
\vspace{0.15cm} 

(1) $Z$ is smoothable by a flat deformation;  

(2) $Y$ is non-singular.} 
\vspace{0.2cm} 

The following conjecture would be true if any projective 
symplectic variety had a {\bf Q}-factorial terminalization 
(cf. Remark 2, (1))  
\vspace{0.2cm}

{\bf Conjecture}. {\em A projective symplectic variety is 
smoothable by a deformation if and only if it admits a 
crepant (symplectic) resolution.}  
\vspace{0.2cm}

{\bf Remark 2}. 
(1) The existence of a {\bf Q}-factorial 
terminalization follows from the ``minimal model conjecture". 
In fact, let $\tilde{Z} \to Z$ be a projective resolution. 
If the minimal model conjecture holds, one can get a {\bf Q}-factorial 
terminalization $Y \to Z$ after running the minimal model 
program for $\tilde{Z}$ over the fixed base $Z$. The 
{\bf Q}-factorial terminalizations are not unique; however, 
they have the same kind of singularities by virtue of 
the proof of Corollary 1. The minimal model conjecture is proved  
by Shokurov [Sh] in dimension 4.  
\vspace{0.12cm} 

(2) It would be interesting to know whether our results 
also hold in the case where $Z$ is a {\em convex} symplectic 
variety. Here we say that $Z$ is {\em convex} 
if there is a birational projective 
morphism from $Z$ to an affine symplectic singularity $X$.     
When $X$ has non-isolated singularity, the usual deformation 
theory of $Z$ (and $X$) has an infinite dimensional nature and 
does not work well. Instead, we should consider a ``Poisson 
deformation" (cf. [G-K]).   
\vspace{0.2cm}

\section{Proof of Main Theorem} 

(i) {\bf Q-factoriality}: 
A normal variety $Z$ is called {\bf Q}-factorial if, for any 
Weil divisor $D$ of $Z$, $mD$ is a Cartier divisor of $Z$ 
for some positive integer $m$.  
Assume that $Z$ is a projective variety with rational singularities. 
Let $\nu : \tilde{Z} \to Z$ be a resolution such that the exceptional 
locus consists of finite number of divisors $E_i$ $(1 \le i \le m)$. 
Denote by $[E_i] \in H^2(\tilde{Z}, {\mathbf Z})$ the corresponding 
classes. By [Ko-Mo, (12.1.6)], $Z$ is {\bf Q}-factorial if and only 
if 
$$ \mathrm{im}[H^2({\tilde Z}, {\mathbf Q}) \to H^0(Z, R^2{\nu}_*{\mathbf Q})] 
 = \mathrm{im}[\oplus{\mathbf Q}[E_i] \to H^0(Z, R^2{\nu}_*{\mathbf Q})].$$ 
We give here a proof of the 
``only if" part. Fix an ample class $\kappa \in H^2(\tilde{Z}, \mathbf{Q})$. 
Let $H^2_0(\tilde{Z}, \mathbf{Q})$ be the primitive part of $H^2(\tilde{Z}, 
\mathbf{Q})$. By definition, we have $$H^2(\tilde{Z}, \mathbf{Q}) 
= H^2_0(\tilde{Z}, \mathbf{Q}) \oplus \mathbf{Q}\kappa.$$ 
Since $Z$ has only rational singularities, the natural map 
$H^2(Z, \mathbf{Q}) \to H^2(\tilde{Z}, \mathbf{Q})$ is an injection 
of mixed Hodge structures ([De]); hence, by the strict compatibility 
of weight filtrations,  
$H^2(Z, \mathbf{Q})$ admits a pure {\bf Q}-Hodge structure of weight 2. 
Then $H^2_0(Z, \mathbf{Q}) := H^2(Z, \mathbf{Q}) \cap 
H^2_0(\tilde{Z}, \mathbf{Q})$ also admits a pure {\bf Q}-Hodge structure 
of weight 2. Note that $H^2_0(Z, \mathbf{C})$ contains $H^{2,0}(\tilde{Z})$ 
and $H^{0,2}(\tilde{Z})$. With respect to the Hodge-Riemann bilinear 
form on $H^2_0(\tilde{Z}, \mathbf{Q})$, we put $V := 
(H^2_0(Z, \mathbf{Q}))^{\perp}$. One can check that 
$$V \cap H^2_0(Z, \mathbf{Q}) = 0,$$ and hence 
$$H^2(\tilde{Z}, \mathbf{Q}) = V \oplus \mathbf{Q}\kappa \oplus 
H^2_0(Z, \mathbf{Q}).$$ 
Note that $H^2_0(Z, \mathbf{Q})$ is mapped to zero by the map 
$H^2(\tilde{Z}, \mathbf{Q}) \to H^0(Z, R^2\nu_*\mathbf{Q})$. 
Since every element of $V$ is of type (1,1), we see that 
$$\mathrm{im}[H^2(\tilde{Z}, \mathbf{Q}) \to 
H^0(Z, R^2\nu_*\mathbf{Q})] = \mathrm{im}[\mathrm{Pic}(\tilde{Z})
\otimes\mathbf{Q} \to H^0(Z, R^2\nu_*\mathbf{Q})].$$ 
Now take a {\bf Q}-Cartier divisor $D$ of $\tilde{Z}$. The push-forward 
$\nu_*D$ is a {\bf Q}-Cartier divisor of $Z$ since $Z$ is {\bf Q}-factorial. 
Now the pull-back $\nu^*(\nu_*D)$ of $\nu_*D$ is well-defined, and 
one has  
$$D = \nu^*(\nu_*D) + \Sigma a_iE_i$$ for some $a_i \in \mathbf{Q}$. 
This implies the ``only if" part.  
\vspace{0.2cm}       

(ii) {\bf Terminal singularity}: 
We recall the definition of a terminal singularity. Assume, for simplicity, 
that $Z$ is a normal variety such that the canonical divisor $K_Z$ is 
a Cartier divisor. 
We say that $Z$ has only terminal (resp. canonical) singularities if, for a resolution 
$\nu: {\tilde Z} \to Z$ such that the exceptional locus consists of 
finite number of divisors $E_i$ $(1 \le i \le m)$, 
$$ K_{{\tilde Z}} = \nu^*K_Z + \Sigma a_i E_i $$ 
with $a_i > 0$ (resp. $a_i \geq 0$) for all $i$. Let $\Sigma \subset Z$ be the 
singular locus. If $Z$ has only terminal singularities, then 
we have $\mathrm{codim}(\Sigma \subset Z) \geq 3$ 
(cf. [Re]). But if, in addition, $Z$ is a symplectic variety, then 
we have $\mathrm{codim}(\Sigma \subset Z) \geq 4$ by [Na 2], [K]. 
Conversely, if $Z$ is a symplectic variety with $\mathrm{codim}(\Sigma 
\subset Z) \geq 4$, then $Z$ has only terminal singularities [Na 2].  
\vspace{0.12cm}

(iii) {\bf Infinitesimal deformations}:  
Assume that $Z$ is a reduced complex space. 
Recall that a deformation $\mathcal{Z}$ of $Z$ over an 
analytic space $S$ with a reference point $0 \in S$ is a 
flat map $\mathcal{Z} \to S$ with a fixed isomorphism 
$\mathcal{Z} \times_{S}\{0\} \cong Z$. 
When $0 \in S$ is a 
spectrum of an Artinian local {\bf C}-algebra with a unique point, 
$\mathcal{Z}$ is called an infinitesimal deformation of $Z$. 
Moreover, an infinitesimal deformation $\mathcal{Z} \to S$ 
is called {\em locally trivial} if, for any Stein open set 
$V \subset Z$, $\mathcal{Z}\vert_V \cong V \times S$. 
  
For $n \geq 0$, we put $A_n := {\mathbf C}[t]/(t^{n+1})$ 
and $S_n : = \mathrm{Spec}A_n$. Suppose that $Z_n \to S_n$ is 
an infinitesimal deformation of $Z$ over $S_n$.  
Define $T_{Z_n/S_n} := 
\underline{\mathrm{Hom}}(\Omega^1_{Z_n/S_n}, \mathcal{O}_{Z_n}).$ 
Let $S_n[\epsilon] := \mathrm{Spec}(A_n[\epsilon]/\epsilon^2)$. 
There is a natural closed immersion $S_n \to S_n[\epsilon]$. 
Moreover, the injection $A_n \to A_n[\epsilon]/\epsilon^2$ 
induces a map $S_n[\epsilon] \to S_n$.  
$\mathrm{Ext}^1(\Omega^1_{Z_n/S_n}, \mathcal{O}_{Z_n})$ 
plays an important role in the deformation theory 
of a {\em reduced} complex space: \vspace{0.12cm}

(iii-a): {\em There is a one-to-one correspondence between 
$\mathrm{Ext}^1(\Omega^1_{Z_n/S_n}, \mathcal{O}_{Z_n})$ 
and the 
set of equivalence classes of infinitesimal deformations of 
$Z_n$ over $S_n[\epsilon]$}. 
\vspace{0.15cm} 

Here two infinitesimal deformations 
$\mathcal{Z}_n$ and $\mathcal{Z}'_n$ of $Z_n$ over $S_n[\epsilon]$ 
are equivalent if there is an $S_n[\epsilon]$-isomorphism between them 
which induces, over $S_n$, the identity map of $Z_n$. 
   
The correspondence is given in the following manner. 
Assume an infinitesimal deformation $\mathcal{Z}_n$ is given. 
By the map $S_n[\epsilon] \to S_n$, we regard $\mathcal{Z}_n$ as 
an $S_n$-analytic space. 
There is an exact sequence 
$$ \epsilon \mathcal{O}_{Z_n} \stackrel{d}\to \Omega^1_{\mathcal{Z}_n/S_n}
\otimes_{A_n[\epsilon]}A_n \to \Omega^1_{Z_n/S_n} \to 0.$$  
Here $d$ is defined by $d(\epsilon f) := fd\epsilon$ for $\epsilon f 
\in \epsilon\mathcal{O}_{Z_n}$. By using the fact that $Z$ is reduced, 
one can show that this $d$ is an injection; hence we have an element 
of $\mathrm{Ext}^1(\Omega^1_{Z_n/S_n}, \epsilon\mathcal{O}_{Z_n})$. 
Since $\mathcal{Z}_n$ is flat over $S_n[\epsilon]$, $\epsilon
\mathcal{O}_{Z_n} \cong \mathcal{O}_{Z_n}$. 

Conversely, let there be an element of 
$\mathrm{Ext}^1(\Omega^1_{Z_n/S_n}, \mathcal{O}_{Z_n}) 
\cong \mathrm{Ext}^1(\Omega^1_{Z_n/S_n}, 
\epsilon\mathcal{O}_{Z_n})$. This element gives an 
$\mathcal{O}_{Z_n}$-module extension of $\Omega^1_{Z_n/S_n}$ 
by $\epsilon\mathcal{O}_{Z_n}$. We consider this extension as 
an $A_n$-module extension. Then 
by the 
homomorphism 
$\mathrm{Ext}^1_{A_n}(\Omega^1_{Z_n/S_n}, \epsilon \mathcal{O}_{Z_n}) \to 
\mathrm{Ext}^1_{A_n}(\mathcal{O}_{Z_n}, \epsilon \mathcal{O}_{Z_n})$
induced by $d: \epsilon\mathcal{O}_{Z_n}  
\to \Omega^1_{Z_n/S_n}$, we have the commutative diagram of exact 
sequences 

\begin{equation}
\begin{CD}
0 @>>> \epsilon\mathcal{O}_{Z_n} @>>> \cal G @>>> \mathcal{O}_{Z_n} @>>> 0 \\
@. @V{id}VV  @V{D}VV @V{d}VV @. \\
0 @>>> \epsilon\mathcal{O}_{Z_n} @>>> \cal F @>{\phi}>> \Omega^1_{Z_n/S_n} @>>> 0 
\end{CD}
\end{equation}

By definition $\mathcal{G} = \{(f, g) \in 
\mathcal{F} \oplus \mathcal{O}_{Z_n}; \phi(f) = 
dg \}$. We give a ring structure to $\mathcal{G}$ by $(f,g)\times (f', g') :
= (g'f + gf', gg')$. By the inclusion $\epsilon\mathcal{O}_{Z_n} \to 
\mathcal{F}$ we 
regard $\epsilon$ as an element of $\mathcal{F}$. Then $\mathcal{G}$ is 
an $A_n [\epsilon]/\epsilon^2$ algebra by defining $\epsilon . (f,g) := 
(g\epsilon, 0)$ for $(f,g) \in \mathcal{G}$. Now put 
$\mathcal{G} = \mathcal{O}_{\mathcal{Z}_n}$. 
Then $(Z_n, \mathcal{O}_{\mathcal{Z}_n})$ is an infinitesimal 
deformation of $(Z_n, \mathcal{O}_{Z_n})$ over $S_n[\epsilon]$. 
\vspace{0.12cm}

{\bf Remark}. Another approach to prove (iii-a) is to use the 
{\em cotangent complex} and the associated tangent cohomology 
(cf. [L-S], [Pa], [Bi, \S. 5]). Let $L^{\cdot}_{Z_n/S_n}$ be 
the cotangent complex for $Z_n \to S_n$. We put  
$$T^i_{Z_n/S_n} := \underline{\mathrm{Ext}}^i(L^{\cdot}_{Z_n/S_n}, 
\mathcal{O}_{Z_n})$$ 
and $$\mathbf{T}^i_{Z_n/S_n} := \mathbf{Ext}^i(L^{\cdot}_{Z_n/S_n}, 
\mathcal{O}_{Z_n}).$$  By the definition of the cotangent complex, 
one has $T^0_{Z_n/S_n} \cong T_{Z_n/S_n}$. Moreover, if $Z$ is 
reduced, then the argument of 
[Ar 2, Prop. 6.1] shows that $T^1_{Z_n/S_n} \cong 
\underline{\mathrm{Ext}}^1(\Omega^1_{Z_n/S_n}, \mathcal{O}_{Z_n})$. 
Since there is a natural map of complexes $L^{\cdot}_{Z_n/S_n} \to 
\Omega^1_{Z_n/S_n}$, we get the maps 
$$\underline{\mathrm{Ext}}^i := 
\underline{\mathrm{Ext}}^i(\Omega^1_{Z_n/S_n}, \mathcal{O}_{Z_n}) 
\to T^i_{Z_n/S_n}$$ and 
$$\mathrm{Ext}^i := 
\mathrm{Ext}^i(\Omega^1_{Z_n/S_n}, \mathcal{O}_{Z_n}) 
\to \mathbf{T}^i_{Z_n/S_n}.$$ 
Therefore, we have a commutative diagram of local-to-global 
exact sequences: 

\begin{equation}
\begin{CD}
H^1(T_{Z_n/S_n}) @>>> \mathrm{Ext}^1 @>>> 
H^0(\underline{\mathrm{Ext}}^1) 
@>>> H^2(T_{Z_n/S_n}) \\
@VVV  @VVV @VVV @VVV \\
H^1(T^0_{Z_n/S_n}) @>>> \mathbf{T}^1_{Z_n/S_n} @>>> 
H^0(T^1_{Z_n/S_n}) @>>> H^2(T^0_{Z_n/S_n})
\end{CD}
\end{equation}
   
By this diagram we see that $\mathrm{Ext}^1(\Omega^1_{Z_n/S_n}, 
\mathcal{O}_{Z_n}) \cong \mathbf{T}^1_{Z_n/S_n}$ when $Z$ 
is reduced.    
\vspace{0.15cm}

Next let us consider a locally trivial deformation $Z_n \to S_n$ 
of $Z$. Locally trivial deformations are controlled by 
$H^1(Z, T_{Z_n/S_n})$ instead of $\mathrm{Ext}^1$: 
\vspace{0.15cm}

(iii-b): {\em There is a one-to-one correspondence between 
$H^1(Z, T_{Z_n/S_n})$ and the 
set of equivalence classes of locally trivial deformations of 
$Z_n$ over $S_n[\epsilon]$.}        
\vspace{0.15cm}

First let us take a Stein open cover $\{V_i\}_{i \in I}$ 
of $Z_n$. Then each locally trivial deformation $\mathcal{Z}_n$ 
of $Z_n$ is constructed by patching $V_i \times_{S_n}S_n[\epsilon]$ 
together in such a way that, over $S_n$, it gives the original 
$Z_n$. Now the situation is the same as the usual deformation of 
a non-singular variety. By the same argument as [SGA 1], Expose.III, 6, 
we have the above correspondence. 
\vspace{0.12cm}
   
(iv) Assume that $Z$ is a {\bf Q}-factorial projective symplectic 
variety with terminal singularities. By (ii), $\mathrm{codim}(\Sigma 
\subset Z) \geq 4$, where $\Sigma := \mathrm{Sing}(Z)$. We define 
$A_n$ and $S_n$ as in (iii). Suppose that $Z_n \to S_n$ is an 
infinitesimal deformation of $Z$ over $S_n$. Let $U$ be the regular 
part of $Z$. We put $U_n := Z_n\vert_U$. 

{\bf Lemma 1}. {\em The natural map $$H^1(Z, T_{Z_n/S_n}) \to 
\mathrm{Ext}^1(\Omega^1_{Z_n/S_n}, \mathcal{O}_{Z_n})$$ 
is a bijection.}  
\vspace{0.12cm} 

{\em Proof}. It suffices to show that the map 
$$\mathrm{Ext}^1(\Omega^1_{Z_n/S_n}, \mathcal{O}_{Z_n}) \to 
H^0(Z, \underline{\mathrm{Ext}}^1(\Omega^1_{Z_n/S_n}, 
\mathcal{O}_{Z_n}))$$ is the zero map. 
Take a point $p \in \Sigma$. Let $Z(p)$ be a small Stein open 
neighborhood of $p \in Z$ and set $U(p) := Z(p) \cap U$. 
We put $Z_n(p) := Z_n\vert Z(p)$ and $U_n(p) := U_n\vert U(p)$.  
We prove that $$\mathrm{Ext}^1(\Omega^1_{Z_n/S_n}, \mathcal{O}_{Z_n}) 
\to H^0(Z_n(p), \underline{\mathrm{Ext}}^1(\Omega^1_{Z_n(p)/S_n}, 
\mathcal{O}_{Z_n(p)}))$$ is the zero map. 
Since $Z$ is Cohen-Macaulay and 
$\mathrm{codim}(\Sigma \subset Z) \geq 3$, by [Ko-Mo, 12.5.6], one 
has  
$$\mathrm{Ext}^1(\Omega^1_{Z_n/S_n}, 
\mathcal{O}_{Z_n})) \cong H^1(U_n, T_{U_n/S_n})$$ 
and $$H^0(Z_n(p), \underline{\mathrm{Ext}}^1(\Omega^1_{Z_n(p)/S_n}, 
\mathcal{O}_{Z_n(p)})) \cong H^1(U_n(p), T_{U_n(p)/S_n}).$$   
Since the symplectic 2-form $\omega$ extends to $\omega_n \in 
H^0(U_n, \Omega^2_{U_n/S_n})$ by [Na], we can identify $T_{U_n/S_n}$ with 
$\Omega^1_{U_n/S_n}$ by this relative 2-form. Now we only have 
to prove that $$H^1(U_n, \Omega^1_{U_n/S_n}) \to 
H^1(U_n(p), \Omega^1_{U_n(p)/S_n})$$ is the zero map. 
Let us consider the Hodge spectral sequences 
$$ H^j(U_n, \Omega^i_{U_n/S_n}) \Rightarrow H^{i+j}(U, A_n)$$ 
and 
$$ H^j(U_n(p), \Omega^i_{U_n(p)/S_n}) \Rightarrow 
H^{i+j}(U(p), A_n),$$ 
where $A_n$ is the constant sheaf with coefficient $A_n$.  
In the latter spectral sequence, we have a decreasing 
filtration $0 \subset F^2_{loc} \subset F^1_{loc} \subset 
F^0_{loc} = H^2(U(p), A_n)$. By the depth argument,  
$E^{0,1}_1 := H^1(U_n(p), \mathcal{O}_{U_n(p)}) 
= 0$ because $Z$ is Cohen-Macaulay and $\mathrm{codim}(\Sigma 
\subset Z) \geq 3$.  So there is an injection 
$Gr^1_{F_{loc}}(H^2(U(p), A_n) \to
H^1(U_n(p), \Omega^1_{U_n(p)/S_n})$. In particular, we have 
a map $F^1_{loc} \to H^1(U_n(p), \Omega^1_{U_n/S_n})$.  
On the other hand, the first spectral sequence degenerates 
at $E_1$ terms for $i+j = 2$ (cf. [Na 1], Lemma 2.7) because 
$\mathrm{codim}(\Sigma \subset Z) \geq 4$. Hence, for 
the decreasing filtration $F^{\cdot}$ of $H^2(U, A_n)$, 
$Gr^1_F(H^2(U, A_n)) = H^1(U_n, \Omega^1_{U_n/S_n})$. 
In particular, we have a {\em surjection} 
$F^1 \to H^1(U_n, \Omega^1_{U_n/S_n})$. The inclusion 
map $U(p) \to U$ induces a map of the second cohomologies, say 
$\iota(p)$  
$$ \iota(p): H^2(U, A_n) \to H^2(U(p), A_n). $$ 
This map preserves the filtrations $\{F^{\cdot}\}$ 
and $\{F^{\cdot}_{loc}\}$; hence it induces a map 
$F^1 \to F^1_{loc}$. By the commutative diagram 
$$ F^1 \longrightarrow\to  H^1(U_n, \Omega^1_{U_n/S_n})$$ 
$$ \downarrow \hspace{3.2cm} \downarrow $$ 
$$ F^1_{loc} \longrightarrow H^1(U_n(p), \Omega^1_{U_n/S_n})$$ 
we only have to prove that the map $\iota(p)$ is the 
zero map. Take a resolution $\nu: \tilde{Z} \to Z$ such 
that $\nu^{-1}(U) \cong U$. Put $\tilde{Z}(p) := \nu^{-1}(Z(p))$. 
Let us consider the commutative diagram  

$$ H^2(\tilde{Z}, A_n) \rightarrow H^2(U, A_n)$$ 
$$ \downarrow \hspace{1.8cm} \downarrow $$ 
$$ H^2(\tilde{Z}(p), A_n) \to H^2(U(p), A_n).$$ 
Here note that $H^2(\cdot, A_n) = H^2(\cdot, {\mathbf C})
\otimes_{\mathbf C}A_n$.  
Since $Z$ is {\bf Q}-factorial, $\mathrm{im}[H^2(\tilde{Z}, A_n) 
\to H^2(\tilde{Z}(p), A_n)]$ is generated by $[E_i]$'s as an 
$A_n$ module by (i), where $E_i$ are exceptional divisors of 
$\nu$. Therefore, the composite 
$$ H^2(\tilde{Z}, A_n) \to H^2(\tilde{Z}(p), A_n) \to 
H^2(U(p), A_n)$$ is the zero map.     
On the other hand, the map $H^2(\tilde{Z}, A_n) 
\to H^2(U, A_n)$ is a surjection by an argument in  
[Na, Proposition 9](cf. (b), p.143) 
because $\mathrm{codim}(\Sigma \subset Z) \geq 4$. Now $\iota(p)$ 
is the zero map by the commutative diagram. Q.E.D.      
\vspace{0.2cm}

Let us start the proof of Main Theorem. 
Let $A$ be an Artinian local $\mathbf{C}$ algebra and 
put $S := \mathrm{Spec}(A)$. 
Recall that a deformation $\mathcal{Z} \to S$ of $Z$ 
is locally 
trivial if, for any Stein open set $V \subset Z$, 
$\mathcal{Z}\vert_V \cong V \times S$.  
We can introduce a {\em locally 
trivial} deformation 
functor of $Z$ 
$$D_{lt} : (\mathrm{Art})_{\mathbf C} \to (\mathrm{Set})$$ 
as a sub-functor of the usual deformation functor $D$ 
of $Z$. One can check that $D_{lt}$ has a pro-representable 
hull $(R_{lt}, \mathcal{Z}_{R_{lt}})$ in the sense of 
Schlessinger [Sch]. Here $R_{lt}$ is a complete local 
$\mathbf{C}$ algebra and $\mathcal{Z}_{R_{lt}}$ is 
a certain projective system of infinitesimal deformations of 
$Z$. Moreover, by the definition of $D_{lt}$, 
one can construct $R_{lt}$ as a suitable quotient of 
the hull $R$ of the usual deformation functor $D$. 
$D$ has $T^1$-lifting property as noted in [Na] 
(for the relationship between $T^1$-lifting property 
and unobstructedness of deformations, see [Ka], [Na 3]).  Hence, by 
Lemma 1 and (iii-b), we see that $D_{lt}$ also has $T^1$-lifting 
property. This implies that $R$ and $R_{lt}$ are both 
complete regular local ring over $\mathbf{C}$. Since 
their cotangent spaces coincide again by Lemma 1, we conclude 
that $R \cong R_{lt}$. This implies that any flat deformation 
of $Z$ is locally trivial.      

\section{Applications} 

In this section we shall generalize Theorem 2.2 of [Na 1] 
to the case of {\bf Q}-factorial terminalizations (cf. 
Theorem 1) and apply Main Theorem and Theorem 1 to get 
Corollaries 1 and 2.  Propositions 1 
and 2 are preliminaries for the proof of Theorem 1. 
 
Let us recall some generalities on deformation 
theory. Let $Y$ be a compact complex space. Then it 
has been proved by [Gr], [Do], [Pa] and others that, there is a 
Kuranishi space $ 0 \in \mathrm{Def}(Y)$ and 
a semi-universal deformation 
$f: \mathcal{Y} \to \mathrm{Def}(Y)$ of $Y$ with $f^{-1}(0) = Y$. 
Assume that $h: Y \to Z$ is a (proper) map of compact 
complex spaces such that $R^1h_*\mathcal{O}_Y = 0$ 
and $h_*\mathcal{O}_Y = \mathcal{O}_Z$. Let 
$\mathcal{Y} \to \mathrm{Def}(Y)$ and 
$\mathcal{Z} \to \mathrm{Def}(Z)$ be the semi-universal 
deformations. 
Then one can construct two maps $h_* : \mathrm{Def}(Y) 
\to \mathrm{Def}(Z)$ and $H_*: \mathcal{Y} \to \mathcal{Z}$ 
in such a way that the following diagram commutes (cf. [Ko-Mo, 11.4])  

\begin{equation}
\begin{CD}
\mathcal{Y} @>{H_*}>> \mathcal{Z} \\ 
@VVV  @VVV  \\
\mathrm{Def}(Y) @>{h_*}>> \mathrm{Def}(Z) 
\end{CD}
\end{equation}

For the convenience of readers, we shall give a sketch of 
the construction of $H_*$ and $h_*$ according to [Ko-Mo, 11.4]. 
First note that $H_*$ and $h_*$ exist in the formal category 
by [Wa]. We want to construct them in the analytic category. 
Let $G \subset Y \times Z$ be 
the graph of $h$. We shall find a graph $\mathcal{G} \subset 
\mathcal{Y}\times\mathcal{Z}$ in such a way that 
$\mathcal{G} \cap (Y \times Z) = G$. Let us consider the 
component $D$ of the relative Douady space of 
$\mathcal{Y} \times \mathcal{Z}/\mathrm{Def}(Y) \times \mathrm{Def}(Z)$ 
containing $[G]$. The projection morphism 
$D \to \mathrm{Def}(Y)$ has a formal section defined by the formal 
contraction morphism around $0 \in \mathrm{Def}(Y)$. 
Then, by [Ar 1, 1.5], one can find an analytic section.    
Note that, when $h$ is a birational map, $(H_*)_t: Y_t \to Z_{h_*(t)}$ 
is a birational map for $t \in \mathrm{Def}(Y)$. 
\vspace{0.15cm}

{\bf Proposition 1}. {\em Let $Z$ be a projective 
symplectic variety. Assume that $h: Y \to Z$ is a 
{\bf Q}-factorial terminalization, that is, 
$h$ is a birational morphism and $Y$ is a 
{\bf Q}-factorial projective symplectic variety 
with terminal singularities.  
Then there is a deformation of $h$: 
$\mathcal{Y}_{\Delta} \stackrel{\tilde{h}}\to 
\mathcal{Z}_{\Delta} \to \Delta$ over a disc 
such that for $t \in \Delta
\setminus{0}$, ${\tilde h}_t: Y_t \to Z_t$ is 
an isomorphism.} 
\vspace{0.2cm}

{\em Proof}. Denote by $(T,0)$ the Kuranishi space of 
$Y$ and let $p: \mathcal{Y} \to T$ be the semi-universal 
family over $T$. By [Na 1], Proposition (2.5), $T$ is smooth. 
Take a projective resolution $\pi: \tilde
{\mathcal{Y}} \to \mathcal{Y}$ in such a way that $\pi$ is 
an isomorphism outside $\mathrm{Sing}(\mathcal{Y})$. 
For a general $t \in T$, 
$\pi_t: \tilde{Y}_t \to Y_t$ is a resolution. 
Moreover, $\pi_t^{-1}(U_t) \cong U_t$ for the regular 
locus $U_t$ of $Y_t$, since $U_t \cap \mathrm{Sing}(\mathcal{Y}) 
= \emptyset$.  Note that 
$\tilde{Y}_t$ is not, in general, a projective variety, but a 
K\"{a}hler manifold. We shall prove that one can deform $Y_t$ 
further to a variety $Y_{t'}$ so that it 
contains no curves. Since there is a map of Kuranishi 
spaces $h_*: T \to \mathrm{Def}(Z)$ as noted above, 
one has a birational 
morphism $Y_{t'} \to Z_{h_*(t')}$. But, $Y_{t'}$ does not 
have any curve; hence this birational morphism should be 
an isomorphism (cf. [Na 1], claim 3, p. 618). 
Since $Y_t$ is a small deformation of a 
symplectic variety $Y$, $Y_t$ is again symplectic by [Na]. 
A symplectic singularity is a canonical singularity (cf. [Be]). 
A canonical singularity is a rational singularity by [El]. 
Therefore, $Y_t$ has only rational singularities.  
Since $Y_t$ has only rational singularities, 
the natural map $H^2(Y_t, \mathbf{C}) 
\to H^2(\tilde{Y}_t, \mathbf{C})$ is an injection of 
mixed Hodge structures ([De]). Hence, by the strict 
compatibility of weight filtrations,   
$H^2(Y_t, \mathbf{C})$ admits a pure Hodge structure of 
weight 2. Let 
$$H^2(Y_t, \mathbf{C}) = H^{2,0}(Y_t)\oplus H^{1,1}(Y_t) 
\oplus H^{0,2}(Y_t)$$ 
be the Hodge decomposition. 
Recall again that there is an injection of pure Hodge structures 
$H^2(Y_t, \mathbf{C}) \to H^2(\tilde{Y}_t, \mathbf{C})$. 
Put $n := \dim Y$ and consider the perfect pairing 
$$ <\:,\:>: H^{2n-2}(\tilde{Y}_t, \mathbf{C}) \times 
H^2(\tilde{Y}_t, \mathbf{C}) \to 
H^{2n}(\tilde{Y}_t, \mathbf{C}) = \mathbf{C}.$$ 
Set $$ V_t := H^2(Y_t, \mathbf{C})^{\perp} \subset 
H^{2n-2}(\tilde{Y}_t, \mathbf{C}).$$ Here we regard 
$H^2(Y_t, \mathbf{C})$ as a subspace of $H^2(\tilde{Y}_t, 
\mathbf{C})$. 
\vspace{0.12cm} 

{\bf Claim 1}. {\em Every element of 
$V_t$ is of type $(n-1,n-1)$.} 
\vspace{0.2cm} 

{\em Proof}. 
If we put $$V^{n-1,n-1}_t 
:= V_t \cap H^{n-1}(\tilde{Y}_t, \Omega^{n-1}_{\tilde{Y}_t}),$$  
then we see that $V^{n-1,n-1}_t = H^{1,1}(Y_t)^{\perp} 
\subset H^{n-1}(\tilde{Y}_t, \Omega^{n-1}_{\tilde{Y}_t})$ 
for the Serre pairing 
$$ <\:,\:>: H^{n-1}(\tilde{Y}_t, \Omega^{n-1}_{\tilde{Y}_t}) 
\times H^1(\tilde{Y}_t, \Omega^1_{\tilde{Y}_t}) 
\to H^n(\tilde{Y}_t, \Omega^n_{\tilde{Y}_t}) = \mathbf{C}.$$ 
Since $h^{n-1,n-1}(\tilde{Y}_t) = h^{1,1}(\tilde{Y}_t)$, 
$\dim V^{n-1,n-1}_t = h^{1,1}(\tilde{Y}_t) - h^{1,1}(Y_t)$. 
On the other hand, $\dim V_t = b_2(\tilde{Y}_t) - b_2(Y_t)$. 
We shall prove that $\dim V_t = \dim V^{n-1,n-1}_t$. 
In order to do this, it suffices to show that all elements of 
$H^2(\tilde{Y}_t)$ of type (2,0) and of type (0,2) are 
mapped to zero by the map $\phi: H^2(\tilde{Y_t}, \mathbf{C}) 
\to H^0(Y_t, R^2{\pi_t}_*\mathbf{C})$. In fact, since the map 
$H^2(Y_t) \to H^2(\tilde{Y}_t)$ is an injection of pure Hodge 
structures of weight 2, if these are proved, then the injection 
is actually an isomorphism at (2,0) and (0,2) parts. 
By the conjugation, it suffices to show that $\phi(\alpha) = 0$ 
for an element $\alpha \in H^2(\tilde{Y}_t)$ of type (2,0). 
We shall prove that, for every point $x \in Y_t$, $\phi(\alpha)_x 
= 0$ in $(R^2(\pi_t)_*\mathbf{C})_x$. Let $\nu: W \to \tilde{Y}_t$ 
be a projective birational morphism such that $W$ is smooth 
and $D := (\pi_t\circ\nu)^{-1}(x)$ is a simple normal crossing 
divisor of $W$. Put $h := \pi_t\circ\nu$. Since 
$R^1\nu_*\mathbf{C} = 0$, $R^2(\pi_t)_*\mathbf{C}$ injects 
to $R^2h_*\mathbf{C}$; hence we have to check that $\alpha$ 
is sent to zero by the composite $$H^2(\tilde{Y}_t, \mathbf{C}) \to 
H^2(W, \mathbf{C}) \to (R^2h_*\mathbf{C})_x (= H^2(D, \mathbf{C})).$$ 
We call this composite $\psi$. Then $\psi$ preserves Hodge filtrations 
of $H^2(\tilde{Y}_t)$ and $H^2(D)$.  
Hence it induces $Gr^2_F(H^2(\tilde{Y}_t)) \to Gr^2_F(H^2(D)) \subset 
H^2(D)$. 
Since $\alpha \in Gr^2_F(H^2(\tilde{Y}_t))$, $\psi(\alpha) \in 
Gr^2_F(H^2(D))$. But, since $Y_t$ has rational singularities at $x$, 
$Gr^2_F(H^2(D)) = H^0(D, \hat{\Omega}^2_D) = 0$ by [Na 1, Lemma (1.2)], 
where $\hat{\Omega}^2_D := \Omega^2_D/(\mathrm{torsion})$. 
Therefore, $\dim V_t = \dim V^{n-1,n-1}_t$.   
Q.E.D.   
\vspace{0.2cm}

The families $$\tilde{\mathcal{Y}} \to \mathcal{Y} \to T$$ 
induces the Kodaira-Spencer maps 
$$\rho_{\tilde{Y}_t}: T_t(T) \to 
H^1(\tilde{Y}_t, \Theta_{\tilde{Y}_t}),$$ 
$$\rho_{Y_t}: T_t(T) \to \mathrm{Ext}^1(\Omega^1_{Y_t}, 
\mathcal{O}_{Y_t}).$$ 
Since $\mathcal{Y} \to T$ is versal at $t \in T$, 
$\rho_{Y_t}$ is surjective. Define a map 
$$ d: H^1(\tilde{Y}_t, \Theta_{\tilde{Y}_t}) \to 
\mathrm{Ext}^1(\Omega^1_{Y_t}, \mathcal{O}_{Y_t})$$ 
as the {\em differential} of the map 
$\mathrm{Def}(\tilde{Y}_t) \to T$. We then have 
$d\circ\rho_{\tilde{Y}_t} = \rho_{Y_t}$. 
Let $U_t \subset Y_t$ be the regular locus of $Y_t$. 
Then $\mathrm{Ext}^1(\Omega^1_{Y_t}, 
\mathcal{O}_{Y_t}) \cong H^1(U_t, \Theta_{U_t})$.   
By this identification, the map $d$ coincides 
with the restriction map 
$H^1(\tilde{Y}_t, \Theta_{\tilde{Y}_t}) 
\to H^1(U_t, \Theta_{U_t})$. 
A symplectic 2-form on $Y_t$ extends to 
a holomorphic 2-form on $\tilde{Y}_t$. 
By this 2-form, we have a map  
$H^1(\Theta_{\tilde{Y}_t}) \to H^1(\Omega^1_{\tilde{Y}_t})$. 
We put 
$$ T_1 := \mathrm{im}[T_t(T) \to H^1(\tilde{Y}_t, 
\Theta_{\tilde{Y}_t})], $$ 
$$ T_2 := \mathrm{im}[T_t(T) \to H^1(\tilde{Y}_t, 
\Omega^1_{\tilde{Y}_t})].$$ 

{\bf Claim 2}. 
$$\dim T_2 \geq \dim H^1(U_t, \Omega^1_{U_t}).$$ 

{\em Proof}. Let us consider the composite 
$$ T_t(T) \to H^1(\tilde{Y_t}, \Omega^1_{\tilde{Y}_t}) 
\to H^1(U_t, \Omega^1_{U_t}).$$
By the symplectic 2-form, we identify $H^1(U_t, 
\Theta_{U_t})$ with $H^1(U_t, \Omega^1_{U_t})$. 
Then the composite above coincides with 
the Kodaira-Spencer map $\rho_{Y_t}$, which is 
surjective. 
\vspace{0.2cm}

Look at the pairing  
$$(\:,\:):  H^1(\tilde{Y}_t, \Theta_{\tilde{Y}_t}) 
\times H^{n-1}(\tilde{Y}_t, \Omega^{n-1}_{\tilde{Y}_t}) 
\to H^n(\tilde{Y}_t, \Omega^{n-2}_{\tilde{Y}_t})$$  

{\bf Claim 3}. {\em In the pairing above, 
$(T_1, V_t) = 0$.}  
\vspace{0.2cm}

{\em Proof}. For $\zeta \in H^1(\tilde{Y}_t, \Theta_{\tilde{Y}_t})$, 
take the corresponding infinitesimal deformation $\tilde{Y}_{t,\epsilon} 
\to \mathrm{Spec}\mathbf{C}[\epsilon]$. Choose an element 
$l \in H^{2n-2}(\tilde{Y}_t, \mathbf{C})$ of type $(n-1,n-1)$. 
Then $(\zeta, l)$ is an obstruction for  
the class $l$ to remain of type $(n-1, n-1)$ under the infinitesimal 
deformation. Now the claim follows from Claim 1 and the fact 
that $R^2p_*\mathbf{C}_{\mathcal{Y}}$ is a local system on $T$
(cf. [Na]\footnote{In [Na], this has been proved under ceratin 
additional conditions: $h^0(Y_{reg}, \Omega^2_{Y_{reg}}) = 1$, 
$h^1(Y, \mathcal{O}_Y) = 0$. 
However, these conditions are not essential. Here we give a more direct 
proof to this fact by using our Main Theorem. Since $\mathcal{Y} \to T$ is 
proper, 
if we choose $T$ small enough, then there is an open covering $\{\mathcal{V}_i\}_{i \in I}$ of $\mathcal{Y}$ with (analytic) trivializations over $T$:   
$\mathcal{V}_i \cong V_i \times T$. Here 
$V_i := \mathcal{Y} \cap \mathcal{V}_i$.  
Let us take a canonical stratification $Y = \coprod_j Y_j$ of 
$Y$ into locally closed smooth subsets such that $\bar{Y}_j = 
\mathrm{Sing}(\bar{Y}_{j-1})$ and 
$Y_0 = Y_{reg}$. Fix a nowhere vanishing $C^{\infty}$-vector field 
$\zeta$ of $T$. 
One can lift $\zeta$ to a $C^{\infty}$-vector field 
$\zeta_i$ of $V_i \times T$ in 
such a way 
that $\zeta_i$ induces a vector field of each $(Y_j \cap V_i) \times T$ 
for each stratum $Y_j$. By using the partition of unity, one can 
glue together $\{\zeta_i\}$ and get a globally defined vector field 
$\tilde{\zeta}$ of $\mathcal{Y}$. This vector field gives us a 
$C^{\infty}$-trivialization of 
$\mathcal{Y}$ over $T$. Hence $R^ip_*\mathbf{C}$ are constant sheaves.}). 
\vspace{0.2cm}

Since $T_2$ is the image of $T_1$ by the map 
$H^1(\Theta_{\tilde{Y}_t}) \to H^1(\Omega^1_{\tilde{Y}_t})$, 
we have $<T_2, V_t> = 0$ by the commutative diagram 
of the pairing maps: 
$$(\:,\:):  H^1(\tilde{Y}_t, \Theta_{\tilde{Y}_t}) 
\times H^{n-1}(\tilde{Y}_t, \Omega^{n-1}_{\tilde{Y}_t}) 
\to H^n(\tilde{Y}_t, \Omega^{n-2}_{\tilde{Y}_t})$$ 
$$ \downarrow \hspace{5.0cm} \downarrow $$ 
$$ <\:,\:>: H^1(\tilde{Y}_t, \Omega^1_{\tilde{Y}_t}) 
\times H^{n-1}(\tilde{Y}_t, \Omega^{n-1}_{\tilde{Y}_t}) 
\to H^n(\tilde{Y}_t, \Omega^n_{\tilde{Y}_t}).$$  
Note that $\mathrm{codim}(V_t \subset H^{n-1,n-1}(\tilde{Y}_t)) 
= h^{1,1}(Y_t)$. On the other hand, since $h^{1,1}(Y_t) 
= h^1(U_t, \Omega^1_{U_t})$(cf. [Na]\footnote{By the argument of 
the footenote (1) above, we see that $h^2(Y) = h^2(Y_t)$ and 
$h^2(U) = h^2(U_t)$ with $U := Y_{reg}$. Since $Y$ is 
{\bf Q}-factorial, one has $H^2(Y) \cong H^2(U)$; hence 
$H^2(Y_t) \cong H^2(U_t)$. Now we get the result by the 
footnote (1) in p. 21 of the eprint version of [Na] (math.AG/0010114).}),  
we see, by Claim 2, that $T_2 = (V_t)^{\perp}$ with respect to 
the pairing $<\:,\:>$ above.   

Now let $C$ be a connected curve on $\tilde{Y}_t$ such that 
$\pi_t(C)$ is not a point. Since $[C] \notin V_t$, we 
see that $[C]^{\perp} \cap T_2$ is a codimension 1 subspace 
of $T_2$ by the perfectness of $<\:,\:>$. 
In particular, we can find an element $\zeta \in T_1$ 
such that $(\zeta, [C]) \neq 0$. 
The $(n-1,n-1)$ classes of curves on $\tilde{Y}_t$ 
constitute a countable subset of $H^{n-1,n-1}(\tilde{Y}_t)$. 
Hence, if we take a generic deformation $\tilde{Y}_{t'}$ of 
$\tilde{Y}_t$, then any effective 1-cycle on $\tilde{Y}_t$ 
which is not contaied in a fiber of $\pi_t$ cannot deform 
sideways. Since $\tilde{\mathcal{Y}} \to T$ is a Kaehler morphism, 
every irreducible component of the Douady space $D(\tilde{\mathcal{Y}}/T)$ 
parametrizing curves 
is proper over $T$ ([Fu 2]). Therefore, all components  
dominating $T$ parametrize the curves which are 
contracted to points by $\tilde{\mathcal{Y}} \to \mathcal{Y}$. 
By the countability of the Douady space ([Fu 3]),  
we now see that, for a generic deformation $\tilde{Y}_{t'}$ of 
$\tilde{Y}_t$, 
any holomorphic curve $C$ on $\tilde{Y}_{t'}$ 
is $\pi_{t'}$-exceptional for $\pi_{t'}: \tilde{Y}_{t'} 
\to Y_{t'}$, that is, $\pi_{t'}(C)$ is a point. 
Then, the variety $Y_{t'}$ has no holomorphic curves. In fact, if 
there is a curve $C$, then one can find a curve $D$ on 
$\tilde{Y}_{t'}$ such that $\pi_{t'}(D) = C$ by using the 
Chow lemma [Hi], which is a 
contradiction. 
\vspace{0.12cm}

Let $h: Y \to Z$ be a {\bf Q}-factorial terminalization 
of a projective symplectic variety $Z$. We put $\Sigma 
:= \mathrm{Sing}(Z)$. There is a closed subset 
$\Sigma_0$ of $\Sigma$ with $\mathrm{codim}(\Sigma_0 
\subset Z) \geq 4$ such that $Z$ is locally a trivial 
family of a rational double point along $\Sigma \setminus 
\Sigma_0$. In the remainder, we put $U := Z \setminus 
\Sigma_0$ and $V:= Y \setminus \mathrm{Sing}(Y)$. 
Note that $h^{-1}(U) \subset V$. 
We denote by $n$ the dimension of $Z$.  
In this 
situation, we have the following generalization of 
[Na 1], Proposition (2.1). 
\vspace{0.2cm}

{\bf Proposition 2}\footnote{{\bf Q}-factoriality of $Y$ 
is not necessary for this proposition}. {\em There is a commutative diagram 
$$\mathrm{Ext}^1(\Omega^1_Y, \mathcal{O}_Y) 
\to H^1(h^{-1}(U), \Theta_{h^{-1}(U)})$$ 
$$\downarrow \hspace{2.90cm} \downarrow $$ 
$$\mathrm{Ext}^1(\Omega^1_Z, \mathcal{O}_Z) \to 
\mathrm{Ext}^1(\Omega^1_U, 
\mathcal{O}_U),$$ 
where horizontal maps are both isomorphisms.} 
\vspace{0.2cm}

{\em Proof}. In the diagram, each space corresponds 
to the set of first order deformations of $Y$, $h^{-1}(U)$, 
$Z$ or $U$, respectively. A first order deformation 
of $Y$ (resp. $h^{-1}(U)$) induces that of 
$Z$ (resp. $U$)(cf. [Na 1], p.614). Each vertical map 
is nothing but this correspondence. On the other hand, 
the horizontal maps are natural ones induced by the 
restriction. The second horizontal map (at the bottom) is 
an isomorphism by [Ko-Mo, 12.5.6] because 
$\mathrm{codim}(\Sigma_0 \subset Z) \geq 3$ and $Z$ 
is Cohen-Macaulay. 
By the same reason, $\mathrm{Ext}^1(\Omega^1_Y, \mathcal{O}_Y) 
\cong H^1(V, \Theta_V).$ 
We only have to prove that the 
first horizontal map is an isomorphism. Set $F := 
Y\setminus h^{-1}(U)$ and $F^0 := V\setminus h^{-1}(U)$. 
Note that $\dim F \leq n-2$ since $h$ is a crepant partial 
resolution of a symplectic variety ([Na 1], Corollary (1.15); 
see also [Na 2], footnote, p.1).  
Let us consider the exact sequence of cohomology 
with coefficient $\mathbf{C}$
$$ H^2_{F^0}(V) \to H^2(V) 
\to H^2(h^{-1}(U)) \to 
H^3_{F^0}(V).$$ 
Let $j: h^{-1}(U) \to V$ and $i: F^0 \to V$ 
be the inclusion maps. 
Then the dual of this sequence coincides with the exact sequence 
of cohomology with compact support: 
$$H_c^{2n-3}(F^0) \to H_c^{2n-2}(h^{-1}(U)) \to 
H_c^{2n-2}(V) \to H_c^{2n-2}(F^0)$$ 
induced by the exact sequence (cf. [Ha,II, Ex.1.19(c)])  
$$ 0 \to j_!(\mathbf{C}_{h^{-1}(U)}) \to 
\mathbf{C}_V \to i_*(\mathbf{C}_{F^0}) \to 0. $$ 
Here we have $H_c^{2n-3}(F^0) = H_c^{2n-2}(F^0) 
= 0$. In fact, the first terms and the third 
terms of the next exact sequences vanish (cf.[L]) since 
$\dim_{\mathbf{C}}F \leq n-2$ and 
$\dim_{\mathbf{C}}(F\setminus F^0) \leq n-3$:  
$$H^{2n-3}(F\setminus F^0) \to H_c^{2n-2}(F^0) 
\to H^{2n-2}(F)$$  
and $$H^{2n-4}(F\setminus F^0) \to H_c^{2n-3}(F^0) 
\to H^{2n-3}(F).$$    
As a consequence, the map 
$$H^2(V) \to H^2(h^{-1}(U))$$ 
is an isomorphism. Moreover, this is a morphism 
of mixed Hodge structures. Since $Y$ is a symplectic 
variety with terminal singularities, $\mathrm{codim}
(\mathrm{Sing}(Y) \subset Y) \geq 4$; hence 
$H^2(V)$ is equipped with a pure Hodge structure 
of weight 2 with 
$\mathrm{Gr}^1_F(H^2(V)) = H^1(\Omega^1_V)$ 
(cf. e-print version math.AG/0010114 
of [Na], footnote at p.21). 
Let us examine the mixed Hodge structure on 
$H^2(h^{-1}(U))$. Let $\nu: \tilde{Y} \to Y$ 
be a projective birational morphism from a 
non-singular variety $\tilde{Y}$ to $Y$ 
such that $\nu^{-1}(h^{-1}(U)) \cong 
h^{-1}(U)$ and such that $\mathrm{Exc}(h\circ\nu)$ is a 
normal crossing variety of $\tilde{Y}$. We put 
$E' := \tilde{Y}\setminus (h\circ\nu)^{-1}(U)$.  
By [De], 
$\mathrm{Gr}^1_F(H^2(h^{-1}(U))) = 
H^1(\tilde{Y}, \Omega^1_{\tilde{Y}}(\mathrm{log}E'))$. 
Let us consider the exact sequence 
of local cohomology
$$H^1_{E'}(\Omega^1_{\tilde{Y}}(\mathrm{log}E')) 
\to H^1(\Omega^1_{\tilde{Y}}(\mathrm{log}E')) 
\to H^1(h^{-1}(U), \Omega^1_{h^{-1}(U)}) 
\to H^2_{E'}(\Omega^1_{\tilde{Y}}(\mathrm{log}E')).$$ 
The first term and the fourth term both vanish.  
This follows from the same argument as the proof 
of [Na 1], Proposition (2.1): by [Na 1, claim 2, p.616],  
$R^k(h\circ\nu)_*\Omega^{n-1}_{\tilde{Y}}(\mathrm{log}
E')(-E')) = 0$ for $k \geq 2$. Then by taking the 
duals of the first term and the fourth term, we get 
the conclusion since $\mathrm{codim}(\Sigma_0 \subset 
Z) \geq 4$ (cf. [Na 1, p.615]).  As a consequence, we 
see that 
$$\mathrm{Gr}^1_F(H^2(h^{-1}(U)) = H^1(h^{-1}(U), 
\Omega^1_{h^{-1}(U)}).$$ 
Since the natural map 
$$H^1(V, \Omega^1_V) \to H^1(h^{-1}(U), \Omega^1_{h^{-1}(U)})$$ 
coincides with the map 
$$\mathrm{Gr}^1_F(H^2(V)) \to \mathrm{Gr}^1_F(H^2(h^{-1}(U))),$$ 
we conclude that it is an isomorphism. 
Finally, by the symplectic 2-form on $V$, this map 
is identified with 
$$H^1(V, \Theta_V) \to H^1(h^{-1}(U), \Theta_{h^{-1}(U)}).$$ 

{\bf Theorem 1}. {\em Let $h: Y \to Z$ be a 
{\bf Q}-factorial terminalization of a projective 
symplectic variety $Z$. Then the Kuranishi spaces 
$\mathrm{Def}(Y)$ and $\mathrm{Def}(Z)$ are both 
smooth of the same dimension. The natural map (cf. [Na 1]) 
$h_*: \mathrm{Def}(Y) \to \mathrm{Def}(Z)$ is a 
finite covering. Moreover, $Z$ has a flat deformation to 
a (non-projective) symplectic variety $Z_t$ 
which is, at the same time, a small deformation of $Y$.} 
\vspace{0.2cm}

{\em Proof}. When $Y$ is non-singular, it is nothing but 
Theorem (2.2) of [Na 1]. In a general case, we modify the 
proof of [Na 1] in the following way. First, Corollary (1.5) 
of [Na 1] holds true when $\tilde{X}$ is a {\bf Q}-factorial 
terminalization of $X$ (here the notation being the same as 
[Na 1]). Next, Proposition (1.6) of [Na 1] follows from a more 
general fact recently proved in [K]. Finally, we 
have to replace Proposition (2.1) of [Na 1] by Proposition 2  
and replace Claim 3 of [Na 1](p.618) by Proposition 1. 
\vspace{0.2cm} 

{\bf Corollary 1}. {\em Let $h: Y \to Z$ and 
$h': Y' \to Z$ be two {\bf Q}-factorial terminalizations 
of a projective symplectic variety $Z$. Assume that 
$Y$ is non-singular. Then $Y'$ is also non-singular.} 
\vspace{0.2cm}

{\em Proof}. Since $h$ is a symplectic resolution, 
the Kuranishi space $\mathrm{Def}(Z)$ of $Z$ is 
non-singular by [Na 1, Theorem (2.2)]. 
Moreover, a general point of $\mathrm{Def}(Z)$ 
parametrizes a non-singular variety. 
Assume now that $Y'$ is singular. By Proposition 1, 
there is a point $t \in \mathrm{Def}(Z)$ parametrizing 
a suitable deformation $Y'_t$ of $Y'$. 
By Main Theorem, $Y'$ does not have any smoothing via deformation; 
hence $Y'_t$ does not have, too. This contradicts 
the fact that a general point of $\mathrm{Def}(Z)$ 
parametrizes a non-singular variety.  
\vspace{0.2cm} 

{\bf Corollary 2}. {\em Let $Z$ be a projective 
symplectic variety. Assume that $Z$ has a {\bf Q}-factorial 
terminalization 
$h: Y \to Z$. Then the following are equivalent. 
\vspace{0.15cm} 

(1) $Z$ is smoothable by a flat deformation;  

(2) $Y$ is non-singular.} 
\vspace{0.2cm}

{\em Proof}. (1) $\Rightarrow$ (2): Suppose that $Y$ is singular. 
Take a small deformation $Y_t$ of $Y$. Then $Y_t$ 
is not smoothable by 
any flat deformation by Main Theorem. Applying 
Proposition 1 to $h: Y \to Z$, we see that $Z_{h_*(t)} 
(\cong Y_t)$ is also non-smoothable. But, by Theorem 1, 
$\mathrm{Def}(Z)$ is non-singular and hence irreducible. 
This contradicts the assumption (1).  
\vspace{0.12cm}

(2) $\Rightarrow$ (1): This is nothing but Proposition 1.

\quad \\
\quad\\

Yoshinori Namikawa \\
Departement of Mathematics, 
Graduate School of Science, Osaka University, JAPAN \\
namikawa@math.sci.osaka-u.ac.jp

\end{document}